\documentclass[11pt]{amsart}

\usepackage{amssymb}
\usepackage[leqno]{amsmath}
\usepackage{amsfonts}
\usepackage{amsopn}
\usepackage{amstext}
\usepackage{amsthm}

\providecommand{\cal}{\mathcal}

\renewcommand{\frak}{\mathfrak}
\newenvironment{pf}{\begin{proof}}{\end{proof}}

\newcommand{\Aaa}{{\cal{A}}}
\newcommand{\Bee}{{\cal{B}}}
\newcommand{\Cee}{{\cal{C}}}
\newcommand{\Ef}{{\cal{F}}}

\newcommand{\Pee}{{\cal{P}}}

\renewcommand{\k}{\kappa}
\newcommand{\su}{\subseteq}

\newcommand{\lam}{{\lambda}}
\newcommand{\al}{\alpha}

\renewcommand{\phi}{\varphi}
\renewcommand{\rho}{\varrho}

\newcommand{\rest}{\restriction}

\newcommand{\loe}{\leqslant}
\newcommand{\goe}{\geqslant}

\newcommand{\subs}{\subseteq}
\newcommand{\sups}{\supseteq}

\newcommand{\ovr}{\overline}

\newcommand{\cf}{\operatorname{cf}}

\newcommand{\Land}{\;\&\;}

\newcommand{\by}{/_}

\newtheorem{tw}{Theorem}[section]
\newtheorem{pblm}[tw]{Problem}
\newtheorem{wn}[tw]{Corollary}
\newtheorem{lm}[tw]{Lemma}

\newtheorem{claim}[tw]{Claim}
\theoremstyle{definition}

\theoremstyle{remark}

\newcommand{\setof}[2]{\{#1\colon #2\}}

\newcommand{\seqof}[2]{\langle #1\colon #2\rangle}

\newcommand{\sn}[1]{\{#1\}} 
\newcommand{\dn}[2]{\{#1,#2\}} 
\newcommand{\map}[3]{#1\colon #2 \to #3} 
\newcommand{\img}[2]{#1[#2]} 
\newcommand{\inv}[2]{{#1}^{-1}[#2]} 

\newcommand{\mad}{\operatorname{MAD}}

\newcommand{\pcf}{\operatorname{pcf}}

\newcommand{\bd}{{\frak{b}}}

\newcommand{\gotha}{{\frak a}}

\theoremstyle{remark}

\newtheorem{fact}[tw]{Fact}

\title[Singular Madness]{On two problems of Erd\H os and Hechler:
New methods in singular Madness}

\author{Menachem Kojman${}^1$} \address{Department of Mathematics,
Ben-Gurion University of the Negev, Beer-Sheva, Israel}
\email{kojman@math.bgu.ac.il}

\thanks{${}^1$ Research partially
supported by an Israeli Science foundation grant no. 177/01}

\author{Wies{\l}aw Kubi\'s}
\address{Department of Mathematics,
Ben-Gurion University of the Negev, Beer-Sheva, Israel
\\ \textit{and} \\ Institute of Mathematics, University of Silesia,
Katowice, Poland}
\email{kubis@math.bgu.ac.il}

\author{Saharon Shelah${}^2$}
\address{Institute of Mathematics, Hebrew University of
        Jerusalem, Israel \\ \textit{and} \\ Department of Mathematics Rutgers
        University, New-Brunswick} \email{shelah@math.huji.ac.il}

\thanks{${}^2$ This research was supported by The
        Israel Science Foundation.  Publication 793.}
 \subjclass[2000]{Primary: 03E10, 03E04, 03E17, 03E35;
 Secondary: 03E55, 03E50}
\keywords{almost disjoint family, singular cardinal, bounding
number, smooth pcf scales}
\begin{document}

\begin{abstract}For an infinite cardinal $\mu$, $\mad(\mu)$
denotes the set of all cardinalities of \emph{nontrivial maximal
almost disjoint families} over $\mu$.

 Erd\H os and Hechler proved in \cite{EH} the consistency of
$\mu\in \mad(\mu)$  for a singular cardinal $\mu$ and asked if it
was ever possible for a singular $\mu$ that $\mu\notin
\mad(\mu)$, and also whether $2^{\cf\mu}<\mu \Longrightarrow
\mu\in \mad(\mu)$ for every singular cardinal $\mu$.

We introduce a new method for controlling $\mad (\mu)$ for a
singular $\mu$ and, among other new results about the structure of
$\mad(\mu)$ for singular $\mu$,  settle both problems
affirmatively.
\end{abstract}

\maketitle

\section{Introduction}

 \subsection{Background} Let $\mu$ be an infinite
cardinal. A family of sets $\Aaa$ is {\em $\mu$-almost disjoint}
({\em $\mu$-ad} for short) if $|A|=\mu=|\bigcup\Aaa|$ for every
$A\in\Aaa$ and $|A\cap B|<\mu$ for every distinct $A,B\in\Aaa$.
$\Aaa$ is {\em maximal $\mu$-almost disjoint} ({\em $\mu$-mad})
if there is no $C\subs\bigcup\Aaa$ such that $\Aaa\cup\{C\}$ is
$\mu$-almost disjoint; in this case we also say that $\Aaa$ is
{\em mad in} $\mu$. It is clear that every $\mu$-almost disjoint
family consisting of fewer than $\cf\mu$ sets is mad in $\mu$;
such a family will be called {\em trivial}. We denote by
$\mad(\mu)$ the set of all cardinalities of nontrivial mad
families in $\mu$. A standard diagonalization argument shows that
$\cf\mu\notin\mad(\mu)$, therefore $\mad(\mu)$ is contained in
the interval of cardinals $[\cf(\mu)^+,2^{\mu}]$.

W. W. Comfort asked (see \cite{EH}) under what conditions it holds
that $\mu\in\mad(\mu)$ for a singular cardinal $\mu$.  P. Erd\H
os and S. Hechler \cite{EH} proved that $\mu\in\mad(\mu)$ if
$\lam^{\cf\mu}<\mu$ for every $\lam<\mu$. Thus, if
$2^{\aleph_0}<\aleph_\omega$ then the interval
$[2^{\aleph_0},\aleph_\omega]$ of cardinals is contained in
$\mad(\aleph_\omega)$.

Erd\H os and Hechler asked in \cite{EH} whether it is consistent
that $\mu\notin\mad(\mu)$ for some singular cardinal $\mu$ and,
more concretely, whether Martin's axiom together with
$2^{\aleph_0}>\aleph_\omega$ implies that
$\aleph_\omega\notin\mad(\aleph_\omega)$. They also asked whether
$2^{\cf\mu}<\mu$ implies $\mu\in\mad(\mu)$ for singular cardinals
$\mu$ other than $\aleph_{\omega}$.

Both problems are settled affirmatively by the general results
below on $\mad(\mu)$ for a singular $\mu$.

\subsection{Notation} Let $\gotha_\mu=\min\mad(\mu)$ and let
$\gotha=\gotha_{\aleph_0}$. For a singular $\mu$ it holds that
$\mad(\cf\mu)\su \mad(\mu)$, therefore $\gotha_\mu\loe
\gotha_{\cf\mu}$.

 A crucial role  in the results is played by two
\emph{bounding numbers}: $\bd_\mu$ and $\bd_{\cf\mu}$.

For every quasi-ordering $(P,\loe)$ with no maximum, the
\emph{bounding number} $\bd(P,\loe)$ is the least cardinality of a
subset of $P$ with no upper bound. For a regular cardinal $\k$,
let $\bd_\kappa$ denote the bounding number of
$(\kappa^\kappa,\loe^*)$, where $f\loe^*g$ means that
$|\setof{i<\kappa}{f(i)> g(i)}|<\kappa$;
 let $\bd=\bd_{\aleph_0}$. It is well known that
$\k<\bd_\kappa\loe\gotha_\kappa$ for a regular cardinal $\kappa$
(for $\kappa=\aleph_0$ see \cite{vD}; the general case is similar)
 and that under Martin's axiom $\bd=2^{\aleph_0}$.

Suppose that $\mu$ is a singular cardinal of cofinality $\k$ and that
$\seqof{\mu_i}{i<\k}$ is a strictly increasing sequence of regular
cardinals with supremum $\mu$. Standard diagonalization shows that
$\bd(\prod_{i<\k}\mu_i,\loe^*)>\mu$. Denote by $\bd_\mu$ the
\emph{supremum of $\bd(\prod \mu_i,\loe ^*)$ over all strictly
increasing sequences of regular cardinals $\seqof{\mu_i}{i<\k}$ with
supremum $\mu$}.

Each  of the following three  relations is  consistent with ZFC:
$\bd<\bd_{\aleph_\omega}$, $\bd=\bd_{\aleph_\omega}$ and
$\bd>\bd_{\aleph_\omega}$.

\subsection{The results} We prove that for every singular cardinal $\mu$:

\medskip
\begin{enumerate}

\item[(1)] $\gotha_\mu\goe\min\dn {\bd_\mu}{\bd_{\cf\mu}}$.

\item[(2)] $\gotha_\mu\loe\lambda<\bd_\mu \Longrightarrow\lambda\in\mad(\mu)$.

\end{enumerate}
\medskip

Thus, if $\bd_{\cf\mu}>\mu$ it follows from $(1)$ that
$\gotha_\mu>\mu$, hence $\mu\notin \mad(\mu)$; and if
$\gotha_{\cf\mu}<\mu$ it follows from (2) that $\mu\in
\mad(\mu)$. In particular:
\medskip
\begin{enumerate}
\item[(b)] $MA+2^{\aleph_0}>\aleph_\omega \Longrightarrow
\aleph_\omega\notin\mad(\aleph_\omega).$

\item[(a)] $2^{\cf\mu}<\mu \Longrightarrow \mu\in\mad(\mu)$ for
every singular $\mu$.
\end{enumerate}
which, respectively,  settle in the affirmative  both problems of
Erd\H os and Hechler from  \cite{EH}.

If one assumes the consistency of large cardinals,
$\bd_{\aleph_\omega}$ can be shifted up arbitrarily high below
$\aleph_{\omega_1}$. Following this with a  ccc forcing for
controlling $\bd$ proves the following:

\begin{enumerate}
\item[(3)] For every  regular $\lambda\in
(\aleph_\omega,\aleph_{\omega_1})$ and  regular uncountable
$\theta\loe \lambda^+$  it is consistent that
$$\mad(\aleph_\omega)=[\theta,\lambda^+].$$
\end{enumerate}

 So, e.g. the following are
consistent:
\begin{itemize}
\item $\mad(\aleph_\omega)=\{\aleph_1,\aleph_2,\dots,\aleph_{\omega+\beta+2}=
2^{\aleph_\omega}\}$ for an arbitrary $\beta<\omega_1$.
\item
$\mad(\aleph_\omega)=\{\aleph_{\omega+\beta+2}\}$ for an arbitrary
$\beta<\omega_1$.
\item
$\mad(\aleph_\omega)=[\aleph_{\omega+\alpha+1},\aleph_{\omega+\beta+2}]$
for arbitrary $\alpha\loe \beta<\omega_1$.
\end{itemize}
And so on.

We refer the reader to the comprehensive list of references in D.
Monk's recent \cite{M}, in which maximal almost disjoint families
are viewed as partitions of unity in the Boolean algebra
$\Pee(\mu)\by {[\mu]^{<\mu}}$.

\subsection{Preliminary facts}
We will use the following facts from \cite{EH}:
\begin{enumerate}
\item[(1)] $\mad(\cf\mu)\subs\mad(\mu)$ and \item[(2)] $\mad(\mu)$ is closed under
singular suprema.
\end{enumerate}
The latter fact is stated in \cite{EH} in a
less general form, so we give a proof here:

\begin{lm}\label{paw} Assume that $\lam=\sup_{i<\theta}\lam_i$, where
$\setof{\lam_i}{i<\theta}\subs\mad(\mu)$ and $\theta<\lam$. Then
$\lam\in\mad(\mu)$.
\end{lm}

\begin{pf} We may assume that $\theta\loe\lam_0$. Let $\Aaa$ be a mad
family in $\mu$ with $|\Aaa|=\lam_0$. Write
$\Aaa=\setof{A_i}{i<\lam_0}$ and for each $i<\theta$ choose a mad
family $\Bee_i$ with $\bigcup\Bee_i=A_i$ and $|\Bee_i|=\lam_i$.
Set
$$\Cee=\bigcup_{i<\theta}\Bee_i\cup\setof{A_j}{\theta\loe j<\lam_0}.$$
Then $|\Cee|=\lam$ and $\Cee$ is mad in $\mu$.
\end{pf}

The following fact will also be used in some proofs.

\begin{lm}\label{smok} Let $\kappa=\cf\mu$ and let $\Aaa$ be a
$\mu$-almost disjoint family of size $\kappa$. Then there exists
a mad family $\Aaa'\sups\Aaa$ such that $|\Aaa'|=\gotha_\mu$ and
$\bigcup\Aaa'=\bigcup\Aaa$.
\end{lm}

\begin{pf} Fix a $\mu$-mad family $\Bee$ with
$|\Bee|=\gotha_\mu$. Choose $\Bee_0=\setof{B_i}{i<\kappa}\subseteq
\Bee$. Let $X=\bigcup_{i<\kappa}B_i$ and define $$\Bee'=\setof{B\cap
X}{X\in\Bee\setminus\Bee_0\Land |X\cap B|=\mu}.$$ Let
$\seqof{A_i}{i<\kappa}$ be a one-to-one enumeration of $\Aaa$. Define
a bijection $\map f{\bigcup\Aaa}X$ so that $\img
f{A_i\setminus\bigcup_{j<i}A_j}=B_i\setminus\bigcup_{j<i}B_j$.
Finally, set $\Aaa'=\Aaa\cup\setof{\inv fB}{B\in\Bee'}$. Observe that
$\Aaa'$ is mad and $|\Aaa'|=\gotha_\mu$.
\end{pf}

\section{Inequalities}
\relax From now on, $\mu$ will always denote a singular cardinal whose
cofinality is denoted by $\kappa$.

\subsection{Bounding numbers and madness in singular cardinals}

\begin{tw}\label{wiewior} For every singular cardinal $\mu$,
\begin{equation}\label{minbb} \gotha_\mu\goe\min
\dn{\bd_\mu}{\bd_{\cf\mu}}.
\end{equation}
\end{tw}

\begin{pf}Let $\k=\cf\mu$.
Suppose to the contrary that
$\gotha_\mu<\min\dn{\bd_\mu}{\bd_\kappa}$ and fix a strictly
increasing sequence of regular cardinals
$\seqof{\mu_i}{i<\kappa}$ with supremum $\mu$ such that
$\bd(\prod_{i<\k}\mu_i,\loe^*)>\gotha_\mu$.

Let $\Aaa=\bigl\{{\sn i\times\mu}:{i<\kappa}\bigr\}$. By Lemma
\ref{smok}, there exists a family $\Bee\subs[\kappa\times\mu]^\mu$
such that $\Bee\cup\Aaa$ is mad in $\mu$, $\Bee\cap
\Aaa=\emptyset$ and $|\Bee|=\gotha_\mu$.

For each $B\in\Bee$, define a function $f_B:\k\to \k$ by
$f_B(i)=\min \{j<\k: |B\cap (\{i\}\times \mu)|<\mu_j\}$. This
function is well defined, since $|B\cap (\{i\}\times \mu)|<\mu$
for each $i<\k$.

Since $|\Bee|=\gotha_\mu <\bd_\k$, there exists a function
$f:\k\to \k$ so that $f_B<^* f$ for all $B\in \Bee$. Without loss
of generality we may assume that $f$ is strictly increasing.

For each $B\in \Bee$, for all but boundedly many $i<\k$ it holds
that $\sup \{\alpha<\mu_{f(i)}: (i,\alpha)\in B\}<\mu_{f(i)}$.
Let $g_B(i)$ be defined by
\[
g_B(i)=
\begin{cases}
0 & \text {if }\sup \{\alpha<\mu_{f(i)}: (i,\alpha)\in
B\}=\mu_{f(i)}\cr
 \sup\{\alpha<\mu_{f(i)}: (i,\alpha)\in B\} & \text{otherwise }
\cr
\end{cases}
\]

For each $B\in \Bee$ the function $g_B$ belongs to $\prod
_{i<\k}\mu_{f(i)}$. Since
$$\bd(\prod_{i<\k} \mu_{f(i)},\loe^*~)\goe
\bd(\prod_{i<\k}\mu_i)>\gotha_\mu$$
 we can fix a function $g\in
\prod_{i<\k}\mu_{f(i)}$ so that $g_B<^* g$ for all $B\in \Bee$.

Define

$$C=\bigcup_{i<\kappa}\sn{i}\times[g(i),\mu_{f(i)}).$$

Clearly, $|C|=\mu$. For each $B\in\Bee$ there exists $j_B<\kappa$
such that $g_B(i)<g(i)$ for all $i>j_B$. This implies that
$\sn{i}\times[g(i),\mu_{f(i)})$ is disjoint from $B$ for all
$i>j_B$. Hence $|B\cap C|\loe\mu_{f(j_B)}<\mu$. Clearly,
$|C\cap(\sn{i}\times\mu)|\loe\mu_{f(i)}<\mu$ for all $i<\k$; so
$\Aaa\cup\Bee\cup\sn{C}$ is $\mu$-almost disjoint, contrary to the
maximality of $\Aaa\cup\Bee$.
\end{pf}

A positive answer to the first question of Comfort, Erd\H os and
Hechler follows now as a corollary:

\begin{wn} If Martin's Axiom holds and
$2^{\aleph_0}>\mu>\cf\mu=\aleph_0$ then $\mu\notin\mad(\mu)$.
\end{wn}

\subsection{Between $\gotha_\mu$ and $\bd_\mu$}
In this section we shall show that $\mad(\mu)$ contains   the
interval of cardinals $[\gotha_\mu,\bd_\mu)$ and even
$[\gotha_\mu,\bd_\mu]$ in the case $\bd_\mu$ is a successor of a
regular cardinal.

\begin{tw}\label{between}For every singular cardinal $\mu$ and
every cardinal $\lambda$,
\begin{equation} \label{inbetween}
\gotha_\mu\loe\lambda<\bd_\mu\Longrightarrow \lambda\in \mad(\mu).
\end{equation}
 If $\bd_\mu$ is a successor of a regular cardinal, then
$\gotha_\mu\loe \bd_\mu\implies \bd_\mu\in \mad(\mu)$.
\end{tw}

To prove  the Theorem it suffices, by Lemma
\ref{paw},  to show that every regular $\lam\in[\gotha_\mu,\bd_\mu)$
belongs to $\mad(\mu)$.

 The proof of
this  will now be divided to two cases. First we prove that every
regular   $\gotha_\mu<\lambda< \mu$ belongs to $\mad(\mu)$. The
proof in this case does not require any specialized techniques. Then
we prove the same for regular $\mu<\lambda<\bd_\mu$ and for
$\bd_\mu$ itself when it is the successor of a regular cardinal.
In this case the proof requires some machinery from pcf theory.

Despite of the technical differences between  both proofs, they
are similar, and could, in fact, be combined to a single proof.
Both follow the same scheme of gluing together $\lam$ different
$\mu$-mad families, each of size $\gotha_\mu$, to a single
$\mu$-mad family of size $\lam$. In the case $\lam<\mu$, a simple
presentation of $\mu$ as a disjoint union of $\lam$ parts works;
in the second part we need to rely on smooth pcf scales to get a
presentation of $\mu$ as an \emph{almost increasing} and
\emph{continuous} union of length $\lam$ of sets of size $\mu$.

\subsubsection{The case $\lambda<\mu$}

\begin{lm} Suppose  $\mu>\cf\mu=\kappa$. Then for every regular
cardinal $\lam$,
 $$\gotha_\mu \loe \lam<\mu\implies \lam\in \mad(\mu).$$
\end{lm}

\begin{pf}
Suppose $\lam$ is regular and $\gotha_\mu\loe \lam <\mu$.  Since
$\gotha_\mu>\kappa=\cf \mu$, $\lam>\kappa$.

Fix a strictly increasing sequence of regular cardinals
$\seqof{\mu_i}{i<\kappa}$ such that $\sup_{i<\kappa}\mu_i=\mu$ and
$\lam<\mu_0$. We will work in $\mu\times\lam$ instead of $\mu$.
Let $S=\setof{\delta<\lam}{\cf\delta=\kappa}$. For each
$\delta\in S$ fix a strictly increasing, continuous sequence
$D_\delta=\seqof{\gamma^{\delta}_i}{i<\kappa}$ with limit
$\delta$ such that $\gamma^\delta_0=0$. Define
$$F_j^\delta=\bigcup\setof{\mu\times\sn \beta}{\gamma^\delta_j\loe\beta<\gamma^\delta_{j+1}}.$$
Thus  $\Ef_\delta=\setof{F_j^\delta}{j<\kappa}$ is a disjoint
family of sets, each set of size $\mu$, which covers
$\mu\times\delta$. Let $\Aaa_\delta\subs[\mu\times\delta]^\mu$ be
such that $\Aaa_\delta\cup\Ef_\delta$ is mad in $\mu\times\delta$,
$\Aaa_\delta\cap \Ef_\delta=\emptyset$ and
$|\Aaa_\delta|=\gotha_\mu$.

Define
$$\Bee=\setof{\mu\times\sn\al}{\al<\lam}\cup\bigcup_{\delta\in S}\Aaa_\delta.$$
Then $|\Bee|=\lam$ and $\Bee\subs [\mu\times\lam]^\mu$.
We will show that $\Bee$ is $\mu$-mad.

First, observe that $\Bee$ is almost disjoint: clearly each element of
$\Aaa_\delta$ is almost disjoint from any set of the form $\mu\times\sn\al$,
because if $\al<\delta$ then $\mu\times\sn\al\subs F_j^\delta$ for $j<\kappa$
such that $\gamma^\delta_j\loe\al<\gamma^\delta_{j+1}$. Finally, consider
$A_i\in\Aaa_{\delta_i}$, $i<2$, with $\delta_0<\delta_1$. Then
$A_0\subs\bigcup_{j<j_0}F_j^{\delta_1}$, where $j_0<\kappa$ is such
that $\delta_0<\gamma^{\delta_1}_{j_0}$. Thus $|A_0\cap A_1|<\mu$.

To see that $\Bee$ is mad fix an arbitrary $Z\in[\mu\times\lam]^\mu$. There
exists a sequence $\seqof{\al_i}{i<\kappa}$ in $\lam$ such that
$$|Z\cap(\mu\times\sn{\al_i})|\goe\mu_i.$$

If $|\setof{\al_i}{i<\kappa}|<\kappa$ then
$|Z\cap(\mu\times\sn\al)|=\mu$ for some $\alpha$. So suppose that
$|Z\cap(\mu\times\sn{\al_i})|<\mu$ for every $i<\kappa$. Taking a
subsequence, we may assume that $\seqof{\al_i}{i<\kappa}$ is
strictly increasing. Let $\delta$ be its supremum. By regularity
of $\lam$,  $\delta\in S$ and therefore
$Z\in[\mu\times\delta]^\mu$. Shrinking $Z$ if necessary, assume
that $Z\subs\bigcup_{i<\kappa}\mu\times\sn{\al_i}$. Then $|Z\cap
F^\delta_j|<\mu$ for every $j<\kappa$. Thus, $|Z\cap A|=\mu$ for
some $A\in\Aaa_\delta$. This completes the proof.
\end{pf}

\begin{wn}\label{interval} Let $\mu>\cf\mu=\kappa$. If
$\gotha_\kappa\loe\mu$ then
$[\gotha_\kappa,\mu]\subs\mad(\mu)$.  In particular, if\/
$2^{\kappa}<\mu$ then $\mu\in\mad(\mu)$.
\end{wn}

Corollary \ref{interval} anwers affirmatively the second question
of Erd\H os and Hechler in \cite{EH}.

\subsubsection{The case $\lambda>\mu$}

 A {\em $(\mu,\lam)$-scale} is a sequence $\ovr f =
\seqof{f_\alpha}{\alpha<\lam}\subs\prod_{i<\kappa}\mu_i$ such that
$\seqof{\mu_i}{i<\kappa}$ is a strictly increasing sequence of regular
cardinals with limit $\mu$, and so that $\alpha<\beta<\lam\implies
f_\alpha<^*f_\beta$ and for every $g\in\prod_{i<\kappa}\mu_i$ there is
$\alpha<\lam$ with $g<^*f_\alpha$.  The relation $f<^*g$ means that
the set $\setof{i<\kappa}{f(i)\goe g(i)}$ is bounded in $\kappa$. If a
$(\mu,\lam)$-scale exists, then $\lam$ must be a regular cardinal
$>\mu$. When $\mu$ is fixed, ``$(\mu,\lam)$-scale" will be abbreviated
by ``$\lam$-scale". A $\lam$-scale $\ovr f$ is {\em smooth} if for
every $\delta<\lam$ with $\cf\delta>\kappa$ the sequence $\ovr f\rest
\delta=\seqof{f_\alpha}{\alpha<\delta}$ is cofinal in
$(\prod_{i<\kappa}f_\delta(i),<^*)$. In this case we say that
$f_\delta$ is an {\em exact upper bound} of $\ovr f\rest\delta$). We
will denote by $[f,g)$ the set $\setof{(i,\alpha)}{i<\kappa \wedge
f(i)\loe \alpha<g(i)}$.

The proof in the present case goes through two steps. First, it
is shown that whenever a smooth $(\mu,\lam)$-scale exists and
$\gotha_\mu<\lam$, it holds that $\lam\in \mad(\mu)$. Then it is
shown that for every $\mu<\lam<\bd_\mu$ there is a smooth
$(\mu,\lam)$-scale and that in case $\bd_\mu$ is a successor of a
regular cardinal there is also a smooth $(\mu,\bd_\mu)$-scale.

\begin{lm} \label{smooth}
Assume $\lam>\mu>\cf\mu=\kappa$ and   there exists a smooth
$(\mu,\lam)$-scale. If  $\gotha_\mu\loe \lam$ then
$\lam\in\mad(\mu)$.
\end{lm}

\begin{pf} Suppose there exists  a  smooth $\lam$-scale
$\seqof {g_\xi}{\xi<\lam}\subseteq \prod_{i<\k}\mu_i$.
 Let $S=\setof{\delta<\lam}{\cf\delta=\kappa}$, and for each
 $\delta\in S$ fix a strictly increasing, continuous,  sequence
$\seqof{\gamma^{\delta}_i}{i<\kappa}$ with limit $\delta$ such that
$\gamma^\delta_0=0$ and put
$D_\delta=\setof{\gamma^{\delta}_i}{i<\kappa}$.

By induction on $\xi<\lam$  we construct a smooth
$\lam$-scale $\ovr f=\seqof{f_\xi}{\xi<\lam}\subseteq
\prod_{i<\k}\mu_i$ which satisfies  the following two conditions:
\begin{enumerate}
\item[(1)] If $\delta<\lam$ is a limit and $\cf\delta\loe \kappa$ then
$f_\delta(i)=\sup_{\xi\in D_\delta}f_\xi(i)$.

\item[(2)] For each $\xi<\lam$ the set $[f_\xi,f_{\xi+1})=\setof{(
i,\alpha)}{f_\xi(i)\loe \alpha<f_{\xi+1}(i)}$ has cardinality $\mu$.
\end{enumerate}

By induction on $\xi<\lam$ we define an increasing and continuous
sequence of ordinals $\zeta(\xi)<\lam$ and a $<^*$-increasing
sequence of functions $f_\xi\in
\prod_{i<\k}\mu_i$ so that $f_\xi=g_{\zeta(\xi)}$ for all $\xi<\lam$
\emph{except} when $\xi$ is limit of cofinality $\loe \kappa$. Then
$\ovr f:=\seqof{f_{\xi}}{\xi<\lam}$ will be a smooth
$\lam$-scale as required.

  At a limit stage $\xi$ of cofinality $\loe\kappa$ let $\zeta(\xi)
=\bigcup_{\xi'<\xi}\zeta(\xi')$ and use condition (1) to define
$f_\xi$; at successor $\xi+1$ choose $\zeta(\xi+1)$ so that
$\max\{f_\xi,g_{\zeta(\xi)}\}<^* g_{\zeta(\xi+1)}$ and (2) holds, and
let $f_{\xi+1}=g_{\zeta(\xi+1)}$. Suppose now that $\xi$ is a limit of
cofinality $>\kappa$. By the smoothness of $\ovr g$, and since
$\seqof{g_{\zeta(\xi')}}{\xi'<\xi}$ is $<^*$-increasing, after
defining $\zeta(\xi)=\bigcup_{\xi'<\xi}\zeta(\xi')$ we get that
$g_{\zeta(\xi)}$ is an exact upper bound of
$\seqof{g_{\zeta(\xi')}}{\xi'<\xi}$. But then $g_{\zeta(\xi)}$ is also
an exact upper bound of $\seqof{f_{\zeta(\xi')}}{\xi'<\xi}$, and we
let $f_\xi=g_{\zeta(\xi)}$.

Let $f_\lam$ be defined on $\kappa$  by $f_\lam(i)=\mu_i$.

\begin{claim}\label{nobigk}
Suppose $\delta\loe \lam$ and $A\su[0,f_\delta)$ has cardinality
$\mu$. If $\cf \delta>\k$ there is some $\delta'<\delta$ so that
$|A\cap [0,f_{\delta'})|=\mu$.
\end{claim}

\begin{pf}Find $g< f_\delta$ so that $\sum_{i<\kappa} |A\cap (i\times g(i))|=\mu$. By
smoothness there exists some $\delta'<\delta$ so that $g<^* g_{\delta'}$.
\end{pf}

For every $\xi<\lam$ let $A_\xi=[f_\xi,f_{\xi+1})$ and let
$\Aaa=\setof{A_\xi}{\xi<\lam}$.  Then $\Aaa\su
\Pee\bigl([0,f_\lam)\bigr)$ is $\mu$-almost disjoint and
$|\Aaa|=\lam$.

 For each $\delta\in S$ and $i<\kappa$ let
$F^\delta_i=[f_{\gamma^\delta_i},f_{\gamma^\delta_{i+1}})$. Then
$\Ef_\delta=\setof{F^\delta_i}{i<\kappa}$ is a $\mu$-almost
disjoint family whose union is, by condition (1)  on
$\ovr f$, equal to $[0,f_\delta)$. Fix a $\mu$-ad family
$\Bee_\delta\subs\Pee([0,f_\delta))$  such that
$|\Bee_\delta|=\gotha_\mu$,
$\Bee_\delta\cup\Ef_\delta$ is $\mu$-mad and $\Bee_\delta\cap
\Ef_\delta=\emptyset$ (by Lemma \ref{smok}).

\begin{claim} \label{supstok}
If $\delta\in S$ and $B\in \Bee_\delta$ then for all
$i<\kappa$ it holds that
$|B\cap [0,f_{\gamma^\delta_i})|<\mu$.
\end{claim}

\begin{pf}
 If not so, let $i_0<\kappa$ be the largest so that $|B\cap
 [0,f_{\gamma^\delta_{i_0}})|<\mu$; $i_0$ exists because $D_\delta$ is
 closed.  Now $|B\cap F^\delta_{i_0}|=\mu$ --- a contradiction.
\end{pf}

 Let $\Bee=\bigcup_{\delta\in S}\Bee_\delta$. Then
$|\Bee|=\gotha_\mu\cdot\lam=\lam$ and therefore $|\Aaa\cup
\Bee|=\lam$. We will show now that $\Aaa\cup\Bee$ is $\mu$-mad.

Suppose that $A=A_\xi\in \Aaa$ and $B\in \Bee_\delta$ for some
$\delta\in S$. If $\xi\goe \delta$ then clearly $|A\cap B|<\mu$
and if $\xi<\delta$, there is some $i<\k$ so that $A_\xi \su^*
F^\delta_i$ and $|A\cap B|<\mu$ follows from Claim \ref{supstok}.

If $B_1\in \Bee_{\delta_1}$ and $B_2\in \Bee_{\delta_2}$ with
$\delta_1<\delta_2$ in $S$, then there is some $i<\k$ so that
$f_{\delta_1}<^*f_{\gamma^{\delta_2}_i}$ and Claim \ref{supstok}
gives $|B_1\cap B_2|<\mu$.

This establishes that $\Aaa\cup \Bee$ is $\mu$-ad. To verify
maximality, let $Z\su [0,f_\lam)$ be arbitrary of size $\mu$. By Claim
\ref{nobigk} the first $\xi\loe \lam$ for which $|Z\cap
[0,f_\xi)|=\mu$ is either a successor or of cofinality $\loe
\k$. Cofinality $< \kappa$ is ruled-out by condition (1) on $\ovr
f$. The case $\xi$ successor implies that $|Z\cap
A_\xi|=\mu$. Finally, in the remaining case $\xi=\delta\in S$, there
is some $B\in \Bee_\delta$ so that $|Z\cap B|=\mu$
\end{pf}

Now the proof of Theorem \ref{between} will be completed by the
following Lemma, whose proof is actually found implicitly in
\cite{CA}. We shall sketch a proof here too.

\begin{lm}\label{smoothbelowb}
Suppose $\mu$ is singular and $\mu<\lambda<\bd_\mu$, $\lambda$
regular. Then there is a smooth $(\mu,\lam)$-scale. If  $\bd_\mu$
is a successor of a regular cardinal, there is also a smooth
$(\mu,\bd_\mu)$-scale.
\end{lm}

\begin{proof} Since $\lambda<\bd_\mu$, there exists a product
$\prod_{i<\kappa}\mu_i$, where $\kappa=\cf\mu$, so that
$\bd(\prod_{i<\k} \mu_i,<^*)>\lambda$.

By Claim 1.3 in \cite{CA}  there exists a  $\lambda$-scale
$\overline f=\langle f_\alpha:\alpha<\lambda\rangle$ in some
$\prod_{i<\kappa}\mu'_i$ such that for all regular $\theta\in
(\kappa,\mu)$ every $\alpha<\lambda$ with $\cf\alpha=\theta$
satisfies that $\overline f\rest \alpha$ is \emph{flat}, that is,
is equivalent modulo the bounded ideal on $\kappa$ to a strictly
increasing sequence of ordinal functions on $\kappa$.

By Lemma 15 in \cite{EUB}, \emph{every} $\alpha<\lambda$ with
$\cf\alpha>\kappa$ satisfies that $\overline f\rest \alpha$ has an
exact upper bound. Now it is clear how to replace $\overline f$
by a smooth $\lambda$-scale.

Suppose now that $\bd_\mu=\lam^+$, $\lam=\cf\lam$. By \cite{351},
4.1, the set $S^{\lam^+}_{<\lam}:=\{\alpha:\alpha<\lam^+ \wedge
\cf\alpha<\lam\}$ is a union of $\lam$ sets, each of which
carries a square sequence. Therefore, $S^{\lam^+}_{<\lam}\in
I[\lam]$. By 2.5 in chapter 1 of \cite{CA}, there exists a
$(\mu,\bd_\mu)$-scale in which all points of cofinality $<\mu$
are flat and therefore a smooth $(\mu,\bd_\mu)$-scale.
\end{proof}

 In contrast to the case of singular $\mu$, let us mention the
following result of A. Blass \cite{B}, which generalizes Hechler's
\cite{Hechler}: it is consistent that $\mad(\aleph_0)=C$, for any
prescribed closed set of uncountable cardinals $C$ which satisfies
that $[\aleph_1,\aleph_1+|C|]\subs C$ and $\lam^+\in C$ whenever
$\lam\in C$ has countable cofinality. For example, by Blass' or by
Hechler's results there are universes of set theory in which
$\mad(\aleph_0)=\dn{\aleph_1}{\aleph_{\omega+1}}$. By Corollary
\ref{interval}, in any universe that satisfies this it holds that
$[\aleph_1,\aleph_{\omega+1}]\subs\mad(\aleph_\omega)$.

Recently Brendle \cite{Brendle}, using techniques from \cite{ad},
proved the consistency of $\gotha=\aleph_\omega$. 

\begin{pblm}
Is it consistent that $\gotha_{\aleph_\omega}=\aleph_\omega$?
\end{pblm}

\section{Consistency results on $\mad(\aleph_\omega)$ from large cardinal axioms}

The inequality (1) can be used to control $\mad(\aleph_\omega)$
by first increasing  $\bd_{\aleph_\omega}$ and then increasing
$\bd$. PCF theory implies that whenever the SCH fails at a
singular cardinal $\mu$, it holds that $\bd_\mu>\mu^+$. On the
other hand, $\bd_\mu$ cannot be changed by a ccc forcing.

Before we state the result, let us recall some pcf terminology.

\[\pcf\{\aleph_n:n<\omega\}=\Bigl\{\bd\bigl(\prod_n
\aleph_n,\loe_I\bigr): I\su \mathcal
P(\omega) \text{ is a proper ideal }\Bigr\}
\]

The relation $<_I$ is defined by $f<_I g\Leftrightarrow
\{n:f(n)\goe g(n)\}\in I$.

$\pcf\{\aleph_n:n<\omega\}$ is an interval of regular cardinals and
has a maximum. For every $\lambda\in
\pcf\{\aleph_n:n<\omega\}$ there exists a \emph{pcf generator}
$B_\lambda\su \omega$ so that the following holds: denote by
$J_{<\lambda}$ the ideal which is generated by $\{B_\theta:\theta\in
\pcf \{\aleph_n:n<\omega\}\wedge\theta< \lambda\}$; then

\[\lambda=\bd\bigl(\prod_n\aleph_n,\loe_{J_{<\lambda}}\bigr)\]

Finally, $(\aleph_\omega)^{\aleph_0}=\max\pcf\{\aleph_n:n<\omega\}\times
2^{\aleph_0}$. Therefore, if $\aleph_\omega$ is a strong limit,
$2^{\aleph_\omega}=\max\pcf\{\aleph_n:n<\omega\}$.

\begin{fact}For every $\beta<\omega_1$ it is consistent
 (from large cardinal axioms) that
 $2^{\aleph_\omega}=\bd_\mu=\aleph_{\omega+\beta+1}$.
\end{fact}

\begin{pf}
Let $V$ be any  universe of set theory in which $\aleph_\omega$
is a strong limit cardinal and
$2^{\aleph_\omega}=\max\pcf\{\aleph_n:n\in
\omega\}=\aleph_{\omega+\beta+1}$ \cite{Sh, GM}.

In $V$, the ideal $J_{<\max\pcf\{\aleph_n:n<\omega\}}$ is proper and
is generated by countably many sets, therefore by simple
diagonalization there exists an infinite $B\su \omega$ so that
$J_{<\max\pcf\{\aleph_n:n<\omega\}}\rest B$ is contained in the ideal
of finite subsets of $B$. Since $\bd(\prod_n
\aleph_n,\loe_{J_{<\max\pcf\{\aleph_n:n<\omega\}}}~)=\aleph_{\omega+\beta+1}$,
it follows that $\bd(\prod _{n\in
B}\aleph_n,\loe^*)=\aleph_{\omega+\beta+1}$, hence
$\bd_{\aleph_\omega}=\aleph_{\omega+\beta+1}$.
\end{pf}

\begin{tw}For every $\beta<\omega_1$ and $\alpha\loe \omega+\beta+2$
it is consistent (from large
cardinals) that $2^{\aleph_\omega}=\aleph_{\omega+\beta+2}$ and
$\mad(\aleph_\omega)=[\aleph_\alpha,\aleph_{\omega+\beta+2}]$.
\end{tw}

\begin{pf} Start from a model $V$ in which $2^{\aleph_0}=\aleph_1$,
$\aleph_\omega$ is strong limit and
$2^{\aleph_\omega}=\aleph_{\omega+\beta+2}$. Such a model
 exists by the previous Fact.

 For every regular $\aleph_{\omega}<\lam \loe \aleph_{\omega+\beta+2}$
 there is a smooth $\lambda$-scale by Lemma \ref{smoothbelowb}.
 Consequently, there is also
 a smooth $\aleph_{\omega+\beta+2}$-scale.

Now apply Theorem \ref{between} to finish the proof.
\end{pf}

By Theorem 5.4(b) in \cite{Baum}, after adding many Cohen
subsets to $\omega_1$, $\max$ $\mad(\aleph_\omega)$ does not
increase by much. Therefore it is consistent to have
$\mad(\aleph_\omega)=[\aleph_1,\aleph_{\omega+\beta+2}]$ as above,
and to have $2^{\aleph_\omega}$ arbitrary large.

\subsection*{Acknowledgements}
The second author would like to thank Uri Abraham for fruitful
discussions of some of the proofs in this paper and to thank Isaac
Gorelic for  useful remarks.

\end{document}